\documentclass[twoside,a4paper,10pt]{article}
\usepackage{latexsym,amssymb}
\usepackage{mltex}
\usepackage{graphicx}
\usepackage{amstext}
\usepackage[frenchb]{babel}
\usepackage[matrix,arrow]{xy}
\def \be{\begin{eqnarray*}}
\def \ee{\end{eqnarray*}}
\def \ben{\begin{enumerate}}
\def \een{\end{enumerate}}
\def \beit{\begin{itemize}}
\def \eeit{\end{itemize}}
\def \bui#1#2{\mathrel{\mathop{\kern 0pt#1}\limits^{#2}}}
\def \buil#1#2{\mathrel{\mathop{\kern 0pt#1}\limits_{#2}}}
\def \bfll{\begin{flushleft}}
\def \efll{\end{flushleft}}
\def \bflr{\begin{flushright}}
\def \eflr{\end{flushright}}

\def \lra{\longrightarrow}

\def \wit{\widetilde}

\newcommand{\pa}[1]{\left(#1\right)}

\title{Remarques sur le spectre de l'op\'erateur de Dirac}
\author{Nicolas Ginoux}

\begin{document}
\maketitle

\noindent{\bf R\'esum\'e:} Nous d\'ecrivons un nouvelle famille d'exemples d'hypersurfaces de la sph\`ere satisfaisant le cas d'\'egalit\'e de la majoration extrins\`eque de C. B\"ar de la plus petite valeur propre de l'op\'erateur de Dirac.\\

\begin{center}{\Large Remarks on the spectrum of the Dirac operator}\end{center}
$ $\\
\noindent{\bf Abstract:} We describe a new family of examples of hypersurfaces in the sphere satisfying the limiting-case in C. B\"ar's extrinsic upper bound for the smallest eigenvalue of the Dirac operator.\\

 $ $\\

Soit $(M^m,g)$ une hypersurface riemannienne compacte orient\'ee (de dimension $m$) de l'espace-mod\`ele $(\wit{M}^{m+1}(c),g)$, o\`u $\wit{M}^{m+1}(0):=\mathbb{R}^{m+1}$ (l'espace euclidien), $\wit{M}^{m+1}(1):=S^{m+1}$ (la sph\`ere ronde \`a courbure sectionnelle $1$), et $\wit{M}^{m+1}(-1):=\mathbb{H}^{m+1}$ (l'espace hyperbolique \`a courbure sectionnelle $-1$). Soit $\lambda_1$ la plus petite valeur propre de l'op\'erateur de Dirac sur $M$ associ\'e \`a $g$ et \`a la structure spinorielle induite. Il est en g\'en\'eral impossible de donner une expression explicite pour cette valeur propre. Il a cependant \'et\'e d\'emontr\'e que la courbure moyenne $H$ ainsi que la courbure ambiante la majorent de mani\`ere naturelle; plus pr\'ecis\'ement, C. B\"ar a montr\'e dans \cite{Baer98} que
\begin{eqnarray}
\lambda_1^2&\leq&\frac{m^2}{4\mathrm{Vol}(M)}\int_MH^2 v_g\,\qquad\quad\;\textrm{ si }\wit{M}=\mathbb{R}^{m+1},\\
\lambda_1^2&\leq&\frac{m^2}{4\mathrm{Vol}(M)}\int_M(H^2+1) v_g\,\quad\textrm{ si }\wit{M}=S^{m+1},\end{eqnarray} et l'auteur a prouv\'e dans \cite{GinHeinD2003,Ginthese} que
\begin{equation}\lambda_1^2\leq\frac{m^2}{4}(\buil{\sup}{M}\hspace{-1mm}H^2-1)\,\qquad\qquad\,\textrm{ si }\wit{M}=\mathbb{H}^{m+1}.\end{equation}  Si (1) et (2) (resp. (3)) montrent une forte analogie avec les majorations de R.C. Reilly \cite{Reil77} (resp. de E. Heintze \cite{Hein88}) pour la premi\`ere valeur propre non nulle du laplacien scalaire, on peut se demander jusqu'o\`u va cette ressemblance, et si en particulier l'\'egalit\'e dans (1), (2) ou (3) n'est atteinte que lorsque $M$ est \emph{minimalement} immerg\'ee dans une sph\`ere g\'eod\'esique, comme c'est le cas pour le laplacien \cite{Reil77,Hein88,ElSI92}. Il a \'et\'e d\'emontr\'e dans \cite{Baer98} (resp. dans \cite{GinHeinD2003}) que, si (1) ou (2) (resp. (3)) est une \'egalit\'e, alors la courbure moyenne $H$ est constante, et que l'\'egalit\'e a lieu si $M$ est une sph\`ere g\'eod\'esique; n\'eanmoins, le probl\`eme de savoir si cette hypersurface est \emph{la seule} \`a jouir de cette propri\'et\'e est \`a ce jour demeur\'e ouvert.\\

\noindent Nous montrons dans cette note que l'\'egalit\'e dans (2) lorsque $M$ n'est pas une sph\`ere g\'eod\'esique n'entraine pas la minimalit\'e de $M$ dans $S^{m+1}$. Nous exhibons une famille d'exemples \`a cet ef\mbox{}fet.\\

\noindent Ce travail a \'et\'e ef\mbox{}fectu\'e \`a l'Institut Max-Planck pour les Math\'ematiques dans les Sciences de Leipzig, que l'auteur tient \`a remercier pour son soutien et son hospitalit\'e. C'est aussi un plaisir de remercier Oussama Hijazi pour son soutien et pour sa lecture critique de l'article.

\section{Structures spinorielles sur les sous-vari\'et\'es}

\noindent Pour les pr\'eliminaires sur la g\'eom\'etrie spinorielle, on se reportera par exemple \`a \cite{LM,BHMM,BFGK,Friedlivre}.\\

\noindent Nous commen\c cons par rappeler quelques \'el\'ements de base sur la restriction de structures spinorielles (cf. aussi \cite{Baer98,Morel2001,Ginthese}). Tout d'abord, si $E$ et $F$ sont deux f\mbox{}ibr\'es vectoriels riemanniens sur une vari\'et\'e donn\'ee, la pr\'esence d'une structure spinorielle sur deux des f\mbox{}ibr\'es $E$, $F$, $E\oplus F$, induit une structure spinorielle sur le troisi\`eme (\cite{Mil65} et \cite{LM}, Prop. 2.15 p. 90). Cela a deux cons\'equences importantes dans les exemples que nous allons \'etudier.\\ D'une part, tout produit riemannien de deux vari\'et\'es spinorielles porte une structure spinorielle naturelle appel\'ee \emph{structure spinorielle produit}. On peut m\^eme montrer que toute structure spinorielle sur un produit riemannien est en fait une structure spinorielle produit.\\ D'autre part, toute sous-vari\'et\'e riemannienne \`a f\mbox{}ibr\'e normal \emph{trivial} d'une vari\'et\'e riemannienne spinorielle porte une structure spinorielle induite. Nous insistons ici sur le fait que cette structure spinorielle d\'epend de la trivialisation du f\mbox{}ibr\'e normal choisie \cite{Baer2001}. C'est pourquoi nous pr\'eciserons toujours ce choix.\\ Nous aurons \'egalement besoin de l'observation suivante. Soient $E_1\lra M'$ et $E_2\lra M'$ deux f\mbox{}ibr\'es vectoriels r\'eels triviaux sur une vari\'et\'e $M'$, de rangs respectifs $n_1$ et $n_2$. Soit $(X_1,\ldots,X_{n_1})$ (resp. $(X_{n_1+1},\ldots,X_{n_1+n_2})$) une trivialisation de $E_1$ (resp. de $E_2$) et supposons la vari\'et\'e $M'$ spinorielle. Toute trivialisation d'un f\mbox{}ibr\'e vectoriel munissant celui-ci de la structure spinorielle triviale, $(X_1,\ldots,X_{n_1})$ (resp. $(X_{n_1+1},\ldots,X_{n_1+n_2})$) induit d'apr\`es ce qui pr\'ec\`ede une structure spinorielle sur le f\mbox{}ibr\'e $TM'\oplus E_1$ (resp. sur $TM'\oplus E_2$). De m\^eme, la trivialisation $(X_1,\ldots,X_{n_1+n_2})$ de $E_1\oplus E_2$ induit une structure spinorielle sur $TM' \oplus \pa{E_1\oplus E_2}$. Mais $TM'\oplus\pa{ E_1\oplus E_2}=\pa{TM'\oplus E_1}\oplus E_2$, donc $TM'\oplus\pa{ E_1\oplus E_2}$ porte aussi la structure spinorielle induite par celle de $TM'\oplus E_1$ et par $(X_{n_1+1},\ldots,X_{n_1+n_2})$. La remarque est que ces deux structures \emph{co\"\i ncident}. Cela s'applique en particulier \`a la restriction de structures spinorielles \`a une sous-vari\'et\'e \`a f\mbox{}ibr\'e normal trivial (voir ci-dessous).

\section{Une famille d'exemples d'hypersurfaces non minimales satisfaisant l'\'egalit\'e dans (2)}

\noindent Nous d\'ecrivons maintenant les exemples cit\'es en introduction. Soient $p$, $q$ deux entiers naturels non nuls, et $r$, $r'$ deux r\'eels strictement positifs. Consid\'erons la vari\'et\'e $M:=S^p(r)\times S^q(r')$, o\`u $S^k(R)$ d\'esigne la sph\`ere ronde de dimension $k$ et de rayon $R$. Notons $g$ la m\'etrique produit sur $M$, et munissons $M$ de la structure spinorielle produit. Rappelons que, lorsque $k>1$, il n'existe qu'une seule structure spinorielle sur $S^k(R)$, alors que si $k=1$, la vari\'et\'e $S^k(R)$ porte deux structures spinorielles non \'equivalentes, une triviale et une non triviale. Cette derni\`ere est celle qui est induite par le plongement standard $S^1(R)\subset\mathbb{R}^2$. Nous choisirons toujours lorsque $k=1$ la structure spinorielle \emph{non triviale} sur $S^k(R)$.\\

\noindent Pour la structure spinorielle produit que l'on a f\mbox{}ix\'ee sur $M$, le spectre du carr\'e de l'op\'erateur de Dirac $D_M$ de $(M,g)$ est l'ensemble des $\mu_k^2+\mu_l'^2$, o\`u $\mu_k$ (resp. $\mu_l'$) d\'ecrit le spectre de l'op\'erateur de Dirac sur $S^p(r)$ (resp. sur $S^q(r')$) (cela est une propri\'et\'e g\'en\'erale des produits riemanniens, cf. \cite{Atiy84}). En utilisant \cite{Sulthese}, le spectre de $D_M^2$ est donc \[\Big\{\frac{1}{r^2}(\frac{p}{2}+n)^2+\frac{1}{r'^2}(\frac{q}{2}+n')^2,\quad n,n'\in\mathbb{N}\Big\}.\] En particulier, la plus petite valeur propre de $D_M^2$ est $\lambda_1^2=\frac{p^2}{4r^2}+\frac{q^2}{4r'^2}$.\\

\noindent Pour un r\'eel $u$ dans $]0,\frac{\pi}{2}[$, posons $r:=\cos(u)$, $r':=\sin(u)$, et renotons $M_u:=S^p(\cos(u))\times S^q(\sin(u))$. La vari\'et\'e $(M_u,g)$ se plonge alors canoniquement et isom\'etriquement dans $S^{p+q+1}(1)$. De plus, la structure spinorielle induite par ce plongement co\"\i ncide avec la structure spinorielle produit que l'on a f\mbox{}ix\'ee sur $M_u$. Pour montrer cela, distinguons deux cas.\\ 

\noindent a) Si $p>1$ et $q>1$, la structure spinorielle de $M_u$ est unique car $M_u$ est simplement connexe.\\

\noindent b) Si $p=1$ (le raisonnement sera identique si $q=1$), la situation est un peu plus compliqu\'ee puisque $M_u$ porte plusieurs structures spinorielles. Consid\'erons le diagramme commutatif suivant dans lequel les fl\`eches sont les inclusions standard:\\

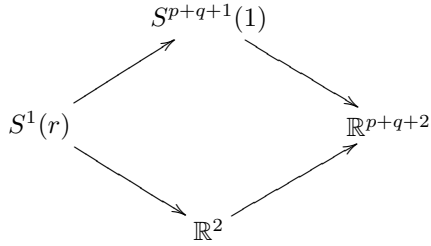
$$\xymatrix{ &S^{p+q+1}(1)\ar[dr]& \\
S^1(r)\ar[ur]\ar[dr]& & \mathbb{R}^{p+q+2}\\
 & \mathbb{R}^2\ar[ur] &}$$
\noindent Nous devons d\'eterminer si la structure spinorielle induite par $S^{p+q+1}(1)$ sur $M_u$ se restreint \`a son tour sur $S^1(r)$ en une structure spinorielle triviale ou non, cette derni\`ere restriction \'etant r\'ealis\'ee au moyen de la trivialisation naturelle (compatible avec les orientations) du f\mbox{}ibr\'e normal de $S^1(r)$ dans $M_u$. Par d\'ef\mbox{}inition, la structure spinorielle induite par $S^{p+q+1}(1)$ sur $M_u$ est donn\'ee par la trivialisation du f\mbox{}ibr\'e normal de $M_u$ dans $S^{p+q+1}(1)$; par suite (cf. paragraphe 1), celle induite par $M_u$ sur $S^1(r)$ est \'egalement induite par $S^{p+q+1}(1)$ sur $S^1(r)$ au moyen de la trivialisation de son f\mbox{}ibr\'e normal dans $S^{p+q+1}(1)$. Or la structure spinorielle de $S^{p+q+1}(1)$ provient elle-m\^eme de la structure spinorielle de $\mathbb{R}^{p+q+2}$ par l'inclusion canonique $S^{p+q+1}(1)\subset\mathbb{R}^{p+q+2}$; la structure spinorielle induite par $M_u$ sur $S^1(r)$ est donc aussi celle induite sur $S^1(r)$ par l'inclusion canonique $S^1(r)\subset\mathbb{R}^{p+q+2}$ au moyen de la trivialisation naturelle du f\mbox{}ibr\'e normal de ce plongement.\\ Mais cette derni\`ere trivialisation est \'egalement obtenue en juxtaposant les trivialisations des f\mbox{}ibr\'es normaux respectifs de $S^1(r)$ dans $\mathbb{R}^2$ et de $\mathbb{R}^2$ dans $\mathbb{R}^{p+q+2}$ (o\`u $\mathbb{R}^2=\mathbb{R}^2\times\{0\}\subset\mathbb{R}^{p+q+2}$); comme $\mathbb{R}^2$ porte \'evidemment la structure spinorielle induite par la trivialisation canonique de son f\mbox{}ibr\'e normal dans $\mathbb{R}^{p+q+2}$, et que le plongement naturel $S^1(r)\subset\mathbb{R}^2$ induit une structure spinorielle \emph{non triviale} sur $S^1(r)$, la structure spinorielle induite sur $M_u$ induit donc \`a son tour une structure spinorielle non triviale sur chacun des facteurs de dimension 1.\\
Par cons\'equent, la premi\`ere valeur propre $\lambda_1$ de l'op\'erateur de Dirac de $M_u$ pour la structure spinorielle induite satisfait \[\lambda_1^2=\frac{p^2}{4\cos(u)^2}+\frac{q^2}{4\sin(u)^2}.\] Or un calcul \'el\'ementaire du carr\'e de la courbure moyenne de $M_u$ dans $S^{p+q+1}(1)$ donne \[H^2=\frac{1}{(p+q)^2}\pa{p\tan(u)-\frac{q}{\tan(u)}}^2,\] d'o\`u 
\[\lambda_1^2=\frac{(p+q)^2}{4}\pa{H^2+1},\] ce qui n'est rien d'autre que l'\'egalit\'e dans (2).\\

\noindent Nous obtenons donc une situation radicalement dif\mbox{}f\'erente de celle du laplacien scalaire, puisqu'il n'existe qu'un seul $u$ dans $]0,\frac{\pi}{2}[$ pour lequel la premi\`ere valeur propre non nulle du laplacien scalaire de $(M_u,g)$ satisfait l'\'egalit\'e dans la majoration de R.C. Reilly \cite{Reil77}, \`a savoir celui pour lequel $M_u$ est \emph{minimale} dans $S^{p+q+1}(1)$. Cette dissemblance, qui vient s'ajouter \`a celle d\'ej\`a observ\'ee par l'auteur \cite{GinHeinD2003} sur la question d'une possible am\'elioration de (3) en une majoration $L^2$ optimale (cf. aussi \cite{Gin2002}), semblerait indiquer une plus grande diversit\'e de r\'esultats pour l'op\'erateur de Dirac dans les probl\`emes d'estimations extrins\`eques des petites valeurs propres. Ce dernier point rejoint par ailleurs la dif\mbox{}f\'erence observ\'ee pour l'\'etude du cas d'\'egalit\'e dans la minoration (intrins\`eque) de T. Friedrich \cite{Fried80} de la plus petite valeur propre de l'op\'erateur de Dirac en fonction de la courbure scalaire (voir \cite{BHMM,Friedlivre} pour ces questions).

\providecommand{\bysame}{\leavevmode\hbox to3em{\hrulefill}\thinspace}

$ $\\

$ $\\

\noindent Nicolas GINOUX\\
Universit\"at Hamburg, FB Mathematik\,-\,SPAD, Bundesstra\ss{}e 55 D-20146 Hamburg\\
fmcw191@math.uni-hamburg.de
\end{document}